\pdfoutput=1
\RequirePackage{ifpdf}
\ifpdf 
\documentclass[pdftex]{sigma}
\else
\documentclass{sigma}
\fi

\usepackage{bbm}
\usepackage{stmaryrd}

\input xy
\xyoption{all}

\newcommand{\Dleft}{[\hspace{-1.5pt}[}
\newcommand{\Dright}{]\hspace{-1.5pt}]}
\newcommand{\SN}[1]{\Dleft #1 \Dright}
\newcommand{\Id}{\mathbbmss{1}}

\newcommand{\rmd}{\textnormal{d}}

\newcommand{\rmh}{\textnormal{h}}

\newcommand{\rml}{\textnormal{l}}
\newcommand{\rmA}{\textnormal{A}}

\newcommand{\rmp}{\textnormal{p}}
\newcommand{\rmb}{\textnormal{b}}
\newcommand{\Lie}{\textnormal{Lie}}

\DeclareMathOperator{\GL}{GL}

\newcommand{\catname}[1]{\textnormal{\texttt{#1}}}

\font\black=cmbx10 \font\sblack=cmbx7 \font\ssblack=cmbx5 \font\blackital=cmmib10 \skewchar\blackital='177
\font\sblackital=cmmib7 \skewchar\sblackital='177 \font\ssblackital=cmmib5 \skewchar\ssblackital='177
\font\sanss=cmss10 \font\ssanss=cmss8 
\font\sssanss=cmss8 scaled 600 \font\blackboard=msbm10 \font\sblackboard=msbm7 \font\ssblackboard=msbm5
\font\caligr=eusm10 \font\scaligr=eusm7 \font\sscaligr=eusm5  \font\fraktur=eufm10
\font\sfraktur=eufm7 \font\ssfraktur=eufm5 
\font\bsymb=cmsy10 scaled\magstep2
\def\all#1{\setbox0=\hbox{\lower1.5pt\hbox{\bsymb
 \char"38}}\setbox1=\hbox{$_{#1}$} \box0\lower2pt\box1\;}
\def\exi#1{\setbox0=\hbox{\lower1.5pt\hbox{\bsymb \char"39}}
 \setbox1=\hbox{$_{#1}$} \box0\lower2pt\box1\;}

\newfam\bifam
\textfont\bifam=\blackital \scriptfont\bifam=\sblackital \scriptscriptfont\bifam=\ssblackital

\newfam\blfam
\textfont\blfam=\black \scriptfont\blfam=\sblack \scriptscriptfont\blfam=\ssblack

\newfam\bbfam
\textfont\bbfam=\blackboard \scriptfont\bbfam=\sblackboard \scriptscriptfont\bbfam=\ssblackboard

\newfam\ssfam
\textfont\ssfam=\sanss \scriptfont\ssfam=\ssanss \scriptscriptfont\ssfam=\sssanss
\def\sss#1{{\fam\ssfam\relax#1}}

\newfam\clfam
\textfont\clfam=\caligr \scriptfont\clfam=\scaligr \scriptscriptfont\clfam=\sscaligr

\newfam\frfam
\textfont\frfam=\fraktur \scriptfont\frfam=\sfraktur \scriptscriptfont\frfam=\ssfraktur

\def\pmb#1{\setbox0\hbox{${#1}$} \copy0 \kern-\wd0 \kern.2pt \box0}
\def\pmbb#1{\setbox0\hbox{${#1}$} \copy0 \kern-\wd0
 \kern.2pt \copy0 \kern-\wd0 \kern.2pt \box0}
\def\pmbbb#1{\setbox0\hbox{${#1}$} \copy0 \kern-\wd0
 \kern.2pt \copy0 \kern-\wd0 \kern.2pt
 \copy0 \kern-\wd0 \kern.2pt \box0}
\def\pmxb#1{\setbox0\hbox{${#1}$} \copy0 \kern-\wd0
 \kern.2pt \copy0 \kern-\wd0 \kern.2pt
 \copy0 \kern-\wd0 \kern.2pt \copy0 \kern-\wd0 \kern.2pt \box0}
\def\pmxbb#1{\setbox0\hbox{${#1}$} \copy0 \kern-\wd0 \kern.2pt
 \copy0 \kern-\wd0 \kern.2pt
 \copy0 \kern-\wd0 \kern.2pt \copy0 \kern-\wd0 \kern.2pt
 \copy0 \kern-\wd0 \kern.2pt \box0}

\mathchardef\za="710B 
\mathchardef\zb="710C 
\mathchardef\zg="710D 
\mathchardef\zd="710E 
\mathchardef\zve="710F 
\mathchardef\zz="7110 
\mathchardef\zh="7111 
\mathchardef\zvy="7112 
\mathchardef\zi="7113 
\mathchardef\zk="7114 
\mathchardef\zl="7115 
\mathchardef\zm="7116 
\mathchardef\zn="7117 
\mathchardef\zx="7118 
\mathchardef\zp="7119 
\mathchardef\zr="711A 
\mathchardef\zs="711B 
\mathchardef\zt="711C 
\mathchardef\zu="711D 
\mathchardef\zvf="711E 
\mathchardef\zq="711F 
\mathchardef\zc="7120 
\mathchardef\zw="7121 
\mathchardef\ze="7122 
\mathchardef\zy="7123 
\mathchardef\zf="7124 
\mathchardef\zvr="7125 
\mathchardef\zvs="7126 
\mathchardef\zf="7127 
\mathchardef\zG="7000 
\mathchardef\zD="7001 
\mathchardef\zY="7002 
\mathchardef\zL="7003 
\mathchardef\zX="7004 
\mathchardef\zP="7005 
\mathchardef\zS="7006 
\mathchardef\zU="7007 
\mathchardef\zF="7008 
\mathchardef\zW="700A 
\mathchardef\zC="7009 

\newcommand{\ra}{\rightarrow}
\newcommand{\lra}{\longrightarrow}
\def\*{{\textstyle *}}
\newcommand{\R}{{\mathbb R}}

\newcommand{\we}{\wedge}
\newcommand{\s}{{\textstyle *}}
\newcommand{\pa}{\partial}
\newcommand{\ti}{\times}
\newcommand{\A}{{\cal A}}

\newcommand{\Ll}{{\pounds}}
\def\lan{\langle}
\def\ran{\rangle}

\def\wt{\widetilde}
\def\Sec{\sss{Sec}}
\def\sJ{{\sss J}}
\def\sT{{\sss T}}

\def\xd{\operatorname{d}}
\def\dt{\xd_{\sT}}
\def\s*{{\scriptstyle *}}
\def\DO{\operatorname{DO}}
\def\cP{\mathcal{P}}

\def\Rt{{\R^\times}}
\def\cG{\mathcal{G}}
\def\cS{\mathcal{S}}

\numberwithin{equation}{section}

\newtheorem{Theorem}{Theorem}[section]
\newtheorem{Proposition}[Theorem]{Proposition}
 { \theoremstyle{definition}
\newtheorem{Definition}[Theorem]{Definition}
\newtheorem{Example}[Theorem]{Example}
\newtheorem{Remark}[Theorem]{Remark} }

\begin{document}

\allowdisplaybreaks

\newcommand{\arXivNumber}{1507.05405}

\renewcommand{\PaperNumber}{059}

\FirstPageHeading

\ShortArticleName{Remarks on Contact and Jacobi Geometry}

\ArticleName{Remarks on Contact and Jacobi Geometry}

\Author{Andrew James BRUCE~$^\dag$, Katarzyna GRABOWSKA~$^\ddag$ and Janusz GRABOWSKI~$^\S$}

\AuthorNameForHeading{A.J.~Bruce, K.~Grabowska and J.~Grabowski}

\Address{$^\dag$~Mathematics Research Unit, University of Luxembourg, Luxembourg} 
\EmailD{\href{mailto:andrewjamesbruce@googlemail.com}{andrewjamesbruce@googlemail.com}}

\Address{$^\ddag$~Faculty of Physics, University of Warsaw, Poland}
\EmailD{\href{mailto:konieczn@fuw.edu.pl}{konieczn@fuw.edu.pl}}

\Address{$^\S$~Institute of Mathematics, Polish Academy of Sciences, Poland}
\EmailD{\href{mailto:jagrab@impan.pl}{jagrab@impan.pl}}

\ArticleDates{Received January 16, 2017, in f\/inal form July 17, 2017; Published online July 26, 2017}

\Abstract{We present an approach to Jacobi and contact geometry that makes many facts, presented in the literature in an overcomplicated way, much more natural and clear. The key concepts are \emph{Kirillov manifolds} and \emph{linear Kirillov structures}, i.e., homogeneous Poisson manifolds and, respectively, homogeneous linear Poisson manifolds. The dif\/ference with the existing literature is that the homogeneity of the Poisson structure is related to a principal ${\rm GL}(1,{\mathbb R})$-bundle structure on the manifold and \emph{not} just to a vector f\/ield. This allows for working with \emph{Jacobi bundle structures} on nontrivial line bundles and drastically simplif\/ies the picture of Jacobi and contact geometry. Our results easily reduce to various basic theorems of Jacobi and contact geometry when the principal bundle structure is trivial, while giving new insights into the theory.}


\Classification{53D05; 53D10; 53D17; 58E40; 58H05}

\section{Introduction}\label{sec:Intro}
There is extensive literature in dif\/ferential geometry devoted to Jacobi structures and derived concepts in which Jacobi structures are presented as generalising Poisson structures. The aim in this paper is to put some order into the f\/ield, and to convince the reader that the properly understood concept of a Jacobi structure is a specialisation of a Poisson structure and not a~ge\-ne\-ralisation. We present an approach to Jacobi and contact geometry which results in drastic simplif\/ication of many concepts, examples and proofs, gives a completely new insight into the theory, as well giving novel discoveries and observations. The main motivation for writing this paper was our observation that many papers in the subject are unnecessarily complicated, because the authors generally ignore, to dif\/ferent extents, the fact that Jacobi geometry is nothing else but homogeneous Poisson geometry on principal $\GL(1,\R)$-bundles ($\Rt$-bundles, for short).

It is common practice to understand the \emph{Jacobi bracket} as a Lie bracket on an algebra of functions on a (smooth) manifold $M$. However, a quick analysis shows its `module nature', and we see that rank $1$ modules (i.e., line bundles) form the natural and proper framework for such structures. In the trivial case, sections of the bundle $\R\ti M\to M$ are identif\/ied with the algebra $\A=C^\infty(M)$ of functions, and the regular $\A$-module structure on $\A$ looks exactly like the multiplication in $\A$, although morphisms in the category of modules are dif\/ferent from these in the category of rings. This we believe is the root of many misunderstandings in Jacobi geometry.

Everyone who works with Jacobi brackets knows there exists something like the \emph{poissonisation} of a Jacobi structure. However, in most cases this is seen as the \emph{poissonisation trick} and used only as a technical tool for proving particular results formulated in the `intrinsic' terms of Jacobi geometry. In consequence, the true landscape of Jacobi geometry, which is actually a~homogeneous Poisson geometry, is hidden in the fog of the `intrinsic Jacobi language'. In our understanding, poissonisation is not a trick, but a genuine framework for Jacobi geometry, and it necessarily comes together with an additional structure of a $\Rt$-bundle. Moreover, the insistence of working in terms of brackets, be they Poisson or Jacobi, often leads to complicated algebraic considerations in which the geometry is completely obscured. In this paper we will work with the corresponding tensor structures and the brackets themselves, which we view as secondary notions, will play no explicit r\^{o}le beyond initial motivation.

Using Jacobi-type brackets on sections of a line bundle $L\to M$ is nothing else but working with `local Lie algebras' in the sense of Kirillov \cite{Kirillov:1976}, so we will call them \emph{Kirillov brackets}. Moreover, the corresponding poissonisations live not on $M\ti\R$ but on the dual bundle $L^*$ with the removed $0$-section, $(L^*)^\ti=L^*\setminus\{ 0_M\}$, which can be recognised as a principal $\GL(1,\R)$-bundle equipped with a homogeneous Poisson tensor. Such structures we will call \emph{Kirillov manifolds} (\emph{Kirillov structures}). We are working with a contact structure if this tensor is actually symplectic. Note that in general by contact structure we mean \emph{non-coorientable contact structure}, that is we do not have a globally def\/ined contact one-form. Working with non-coorientable contact structures forces one to consider non-trivial line bundles and so one has to leave the setting of standard Jacobi geometry. We generally use the identif\/ication $\GL(1,\R)\simeq\R^\times$, where $\R^\times=\R\setminus\{ 0\}$ is the group of multiplicative reals. The appearance of the non-connected group~$\R^\ti$ is forced by the fact that real line bundles over $M$ are classif\/ied by a~$\mathbb{Z}_2=\R^\ti/\R_+$ cohomology of~$M$.

It is well known that the choice of one of equivalent def\/initions inf\/luences strongly our way of thinking and makes the formulations of some concepts and generalisations easier or harder, depending on the choice made. We will insist in this paper in understanding the `poissonisation' as the genuine Jacobi structure with all consequences of this choice. In this sense, the properly understood concept of a Jacobi structure is a specialisation rather than a~generalisation of a~Poisson structure. However, in contrast to the standard def\/inition, we understand the homogeneity not as associated with a vector f\/ield, but rather a principal bundle structure. Thus, our Poisson structure is `homogeneous' in the sense that it is homogeneous of degree $-1$ on a certain $\Rt$-bundle. The appearance of this principal bundle structure is absolutely fundamental for the whole picture. In other words, the proper playground for Jacobi geometry will be the category of \emph{Kirillov manifolds}, i.e., $\Rt$-bundles equipped with a homogeneous Poisson structure. Moreover, all derived concepts like `Jacobi algebroid', `Jacobi bialgebroid', `Jacobi/contact groupoid', etc., should be understood as the corresponding objects in Poisson geometry, equipped additionally with a principal $\R^\ti$-action which is compatible with the other structures. The only thing to be decided is a reasonable notion of compatibility.

Compatibility with a vector bundle structure can be described, in the spirit of viewing linear structures as def\/ined by a certain action of the monoid $(\R,\cdot)$ of multiplicative reals \cite{Grabowski:2006,Grabowski2012}, as commutation of $\R$- and $\R^\times$-actions. As a side remark, actions of the monoid of complex numbers were explored in \cite{Jozwikowski:2016}. Following the observations made in \cite{Bursztyn:2014} (cf.\ \cite{Bruce:2014c,Bruce:2014b, Bruce:2014}), the compatibility with a Lie algebroid/groupoid structure is described as an $\R^\ti$-action by the Lie algebroid/groupoid morphisms, etc. We will call such Lie groupoids $\mathbb{R}^{\times}$-\emph{groupoids}, or more generally, for any Lie group $G$, we have the notion of a $G$-\emph{groupoid}. The inf\/initesimal objects are \emph{$G$-algebroids}. In particular, $\Rt$-algebroids turn out to be the objects dual to linear Poisson manifolds equipped with a compatible principal $\Rt$-action (\emph{linear Kirillov structures}).

The compatibility of a symplectic or Poisson structure with the $\R^\ti$-action is expressed in terms of homogeneity. We must again stress that this homogeneity is not def\/ined in terms of a~vector f\/ield like in \cite{Crainic:2007,Libermann:1995,Marle:1991}, but in terms of the principal $\R^\ti$-action. The fundamental vector f\/ields def\/ine only the action of the connected component in the group. Another important ingredient of our framework, not really exploited in the literature, is the use of tangent and phase lifts of principal $\Rt$-bundle structures.

\looseness=-1 We stress that our def\/inition of a contact groupoid, i.e., a $\mathbb{R}^{\times}$-groupoid equipped with a~homogeneous and multiplicative symplectic form, turns out to be equivalent to the def\/inition of Dazord~\cite{Dazord:1995}. The corresponding objects are called by many authors \emph{conformal contact groupoids}, for example~\cite{Crainic:2007}. The f\/irst and frequently used def\/inition, presented in~\cite{Kerbrat:1993}, is less general and involves an arbitrary multiplicative function, that is due to the fact that in this approach contact bundles are forced to be trivial. Also the def\/inition proposed by Libermann~\cite{Libermann:1993} requires a~contact form. There are, however, no doubts that allowing for contact structures which do not come from a global contact form is fundamental for the completeness and elegance of the theory. On the other hand, although Dazord's def\/inition is general and simple, it does not have any direct extension to Jacobi groupoids. Crainic and Salazar~\cite{Crainic:2015} also realise that there are artif\/icial complications that arise when insisting on working with trivial line bundles, but their way of solving the problem is completely dif\/ferent to ours and makes use of Spencer operators.

Our framework produces, for any Lie groupoid $\mathcal{G}$, canonical examples of contact groupoids $C(\mathcal{G}) \subset \sT^{*}\mathcal{G}$, consisting of all covectors that vanish on vectors that are tangent to the source or target f\/ibres. Interestingly, such examples are somewhat universal and generic; every contact groupoid can be represented or realised as a contact subgroupoid of a canonical one. We will of course be more precise in due course.

Among the main results of this paper:
\begin{itemize}\itemsep=0pt
\item We give the description of any $G$-groupoid in terms of a `splitting' into a product of the $G$-bundle of units and the reduced groupoid (Theorem \ref{generalsplit}). In the simplest case, this description reduces to the well-known groupoid extension by the additive~$\R$ with a~help of a multiplicative function.
\item We show that $G$-algebroids and $G$-groupoids are related via Lie theory (Theorem \ref{thm:lie theory01}).
\item We observe that the dual objects of linear Kirillov structures are exactly $\Rt$-algebroids (Proposition \ref{prop:Kirillov Lie}).
\item We present the Lie theory of Kirillov manifolds and contact groupoids (Proposition \ref{prop:lie theory2}).
\item We prove that all contact groupoids have a realisation as a contact subgroupoid of a~cano\-ni\-cal contact groupoid (Theorem \ref{thm:Dazord}).
\end{itemize}
To sum up, our intention is to show how the setting of homogeneous Poisson geometry simplif\/ies various concepts and problems of Jacobi and contact geometry. We would like to emphasise that, besides a number of new observations, the novelty of this paper lies mainly in the underlying approach that unif\/ies and generalises various points of views, as well as establishes the proper language for the Jacobi and contact geometry. This results in a drastic simplif\/ication of proofs of various important facts spread over the literature and clarif\/ication of the strategies used. We hope that our work will put a new light on Jacobi and contact geometry as a whole, even if separate observations may seem to be known to the reader.

\begin{Remark} The picture of Jacobi geometry -- or better Kirillov geometry -- in terms of $\Rt$-bundles and homogeneous Poisson structures allows for a natural generalisation thereof to the world of $L_{\infty}$-algebras via replacing the Poisson structure with a higher Poisson structure (also known as a~$P_{\infty}$-structure). This leads to the notions of \emph{higher Kirillov manifolds} and \emph{Kirillov structures up to homotopy}, see~\cite{Bruce:2016}.
\end{Remark}

To emphasise the tight parallels with classical constructions/results from Poisson geometry and to help guide the reader, we present the following table.

\begin{table}[h!]\centering\renewcommand{\arraystretch}{1.25}
\begin{tabular}{ |p{7cm}||p{7cm}| }
 \hline
 \multicolumn{2}{|c|}{Poisson geometry and Kirillov geometry} \\
 \hline
 \textbf{Poisson geometry} & \textbf{Kirillov geometry} \\
 \hline
\emph{Poisson structure}: a bi-vector f\/ield $\zL$ on a smooth manifold $M$, that satisf\/ies $\SN{\zL, \zL} =0$, where the bracket is the Schouten--Nijenhuis bracket on $M$. & \emph{Kirillov structure}: a homogeneous bi-vec\-tor f\/ield $\zL$ on a principal $\R^\ti$-bundle $(P, \rmh)$, that satisf\/ies $\SN{\zL, \zL} =0$, where the bracket is the Schouten--Nijenhuis bracket on $P$.\\
\hline
\emph{Poisson bracket}: Lie bracket on functions on $M$. & \emph{Kirillov bracket}: Lie bracket on sections of the line bundle $L$ associated with $(P, \rmh)$. \\
\hline
\emph{Non-degenerate Poisson structure}: symplectic geometry & \emph{Non-degenerate Kirillov structure}: contact geometry.\\
\hline
\emph{Linear Poisson manifold}: Lie algebroid (via dualisation). & \emph{Linear Kirillov manifold}: $\R^\ti$-algebroid (via dualisation).\\
\hline
\emph{Associated Lie algebroid}: the cotangent bundle of a Poisson manifold. & \emph{Associated $\R^\ti$-algebroid}: the cotangent bundle of a Kirillov manifold with the cotangent lift of the $\R^\ti$-action.\\
\hline
\emph{Poisson groupoid}: Lie groupoid with a~mul\-ti\-plicative Poisson structure. Symp\-lec\-tic groupoids are the non-degenerate case. & \emph{Kirillov groupoid}: $\R^\ti$-groupoid with a~mul\-tiplicative Kirillov structure. Contact groupoids are the non-degenerate case.\\
\hline
\emph{Integrating objects}: symplectic groupoids. & \emph{Integrating objects}: contact groupoids.\\
\hline
\end{tabular}
\end{table}

\textbf{Arrangement of paper.} In Section~\ref{sec:PrincipalRbundles} we recall the equivalence of line bundles and principal $\mathbb{R}^{\times}$-bundles. We also remind the reader of the tangent and phase lift of $\mathbb{R}^{\times}$-actions as these will feature heavily throughout this work. In Section~\ref{sec:Kirillov} we show how $\mathbb{R}^{\times}$-bundles appear in an essential way when dealing with Kirillov brackets and contact structures. We then turn our attention to Lie groupoids that have a compatible action of a Lie group upon them in Section~\ref{sec:Ggroupoids}. In this section we examine the structure of such Lie groupoids and show how $\Rt$-algebroids and $\mathbb{R}^{\times}$-groupoids are related via the Lie functor. In Section~\ref{sec:contact and jacobi} we proceed to the main concept of this work: Kirillov and contact groupoids within the framework of $\mathbb{R}^{\times}$-groupoids.

\section[Principal $\R^\times$-bundles]{Principal $\boldsymbol{\R^\times}$-bundles}\label{sec:PrincipalRbundles}

\subsection[Line bundles and principal $\mathbb{R}^{\times}$-bundles]{Line bundles and principal $\boldsymbol{\mathbb{R}^{\times}}$-bundles}\label{subsec:Lines and R^x}

For a vector bundle $E\to M$, with $E^\ti$ we will denote the total space $E$ with the zero-section removed, $E^\ti=E\setminus\{ 0_M\}$. The latter is no longer a vector bundle, but a principal $\R^\ti$-bundle, $E^\ti\to \cP(E)=E^\ti/\R^\ti$, where the principal $\R^\ti$-action, $\rmh\colon \mathbb{R}^{\times} \times E^{\times} \rightarrow E^{\times}$, depends on multiplication by non-zero reals. The bundle $\cP(E) \rightarrow M$ is known as the \emph{projectivisation} of~$E$.

Since $\R^\ti=\GL(1,\R)$, in the case of vector bundles of rank 1, i.e., for line bundles, we have the following.
\begin{Proposition}
The association $L\mapsto L^\ti$ establishes a one-to-one correspondence between line bundles over $M$ and principal $\mathbb{R}^{\times}$-bundles over $M$.
\end{Proposition}

Denote the converse of the above association with $P\mapsto \bar P$. In other words, $\bar P^\ti=P$. Let us introduce the notation $P^+=\bar P^*$.

The fundamental vector f\/ield $\Delta_{P}$ of the $\mathbb{R}^{\times}$-action on $P = L^{\times}$ is nothing more than the Euler vector f\/ield $\Delta_{L}$ on $L$ restricted to $L^{\times}$. We will refer to $\Delta_{P}$, by some abuse of nomenclature, as the \emph{Euler vector field} of~$P$ and the $\mathbb{R}^{\times}$-action as the \emph{homogeneity structure} on $P$.

By employing $t$ as the standard coordinate on $\mathbb{R}$, and so $\mathbb{R}^{\times}$, we can understand $t$ as the f\/ibre coordinate of $P$ in some f\/ixed local trivialisation. That is, $(x^{a}, t)$ can serve as local coordinates on~$P$, where $(x^{a})$ are understood as local coordinates on $M$; such coordinates we will call \emph{homogeneous coordinates}. With respect to homogeneous coordinates the Euler vector f\/ield on~$P$ is simply $\Delta_{P} = t \partial_{t}$.

\subsection{Tangent and phase lifts}
The base manifold of a principal $\R^\ti$-bundle $P$ will be generally denoted with $P_0$, $P_0=P/\R^\ti$. A~fundamental observation is that
principal $\R^\ti$-actions on~$P$ can be canonically lifted to principal $\R^\ti$-actions on~$\sT P$ (\emph{tangent lifts}) and~$\sT^* P$ (\emph{phase lifts}), see, e.g.,~\cite{Grabowski2013}.
\begin{Proposition}\label{prop:tangent lift} Let $\zp\colon P \rightarrow P_0$ be a principal $\mathbb{R}^{\times}$-bundle with respect to an action~$\rmh$. Then,
\begin{itemize}\itemsep=0pt
\item[$(a)$] $\sT P$ is also canonically a principal $\mathbb{R}^{\times}$-bundle, with respect to the action
\begin{gather*}
(\sT \rmh)_{s} := \sT(\rmh_{s}).
\end{gather*}
The base of the corresponding fibration is the Atiyah bundle $\sT P/\R^\ti$, whose sections are interpreted as $\R^\ti$-invariant vector fields on $P$ or, equivalently, as the bundle $\DO^1(P^+,P^+)$ of first-order linear differential operators on the line bundle $P^+$ with values in $P^+$.
\item[$(b)$]
$\sT^* P$ is also canonically a principal $\mathbb{R}^{\times}$-bundle, with respect to the action
\begin{gather*}
(\sT^* \rmh)_{s} := s\cdot (\sT\rmh_{s^{-1}})^*.
\end{gather*}
The base of the corresponding fibration is the first jet bundle $\sJ^1P^+$ of sections of the line bundle~$P^+$.
\end{itemize}
\end{Proposition}

In homogeneous local coordinates $(t, x^{a})$ on $P$, the naturally induced coordinates on the tangent bundle are $(t, x^{a}, \dot{t}, \dot{x}^{b})$. Then,
\begin{gather*}(\sT \rmh)_{s}\big(t, x^{a}, \dot{t}, \dot{x}^{b}\big) = \big(s t, x^{a}, s \dot{t}, \dot{x}^{b}\big).\end{gather*}
 However, these induced coordinates are not necessarily the most convenient. Note that as $t \neq 0$, we can multiply or divide the coordinate functions by $t$ without any ill ef\/fects. In particular, we can make the change $\dot t \mapsto \dot{\textbf{t}} := t^{-1} \dot t$, and then employ coordinates $(t, x^{a}, \dot{\textbf{t}}, \dot{x}^{b})$ on $\sT P$. The admissible changes of coordinates -- which the reader can quickly verify -- are then of the form
\begin{alignat*}{3}
& x^{a'} = x^{a'}(x), \qquad && t' = \phi(x)t,& \\
& \dot{\textbf{t}}' = \dot{\textbf{t}} + \phi^{-1}(x) \dot{x}^a \frac{\partial \phi(x)}{\partial x^a},\qquad && \dot{x}^{b'} = \dot{x}^{a}\frac{\partial x^{b'}(x) }{\partial x^a}.&
\end{alignat*}
Then, with this convenient class of local coordinates we see that the base of $\sT P$, which is the Atiyah bundle $\big(\sT P\big)_0 = \sT P/\R^\ti$, comes with local coordinates $(x^{a}, \dot{\textbf{t}}, \dot{x}^{b})$.

Sections of the Atiyah algebroid $\sT P/\Rt$ of the principal bundle $(\sT P,\sT h)$ can be identif\/ied with invariant vector f\/ields on $P$; locally of the form $X=\alpha(x)t\pa_t+\beta^a(x)\pa_{x^a}$. If $\zs(x)=f(x)$ is a section of $P^+$, then we can apply $X$ to the homogeneous function $\zi_\zs(x,t)=tf(x)$ on $P$ associated with $\zs$,
\begin{gather*}X(\zi_\zs)(x,t)=t\left(\alpha(x)f(x)+\beta^a(x)\frac{\pa f}{\pa x^a}(x)\right).
\end{gather*}
As $X(\zi_\zs)$ is homogeneous, it is of the form $\zi_{\bar X(\zs)}$ for some section $\bar X(\zs)$ of $P^+$. The map $\zs\mapsto\bar X(\zs)$ is a f\/irst-order dif\/ferential operator on $P^+$ with values in $P^+$ of the local form $\bar X=\beta^a(x)\pa_{x^a}+\alpha(x)$.

Similarly, for the naturally induced coordinates $(t,x^a,p,p_b)$ on the cotangent bundle we have
\begin{gather*} (\sT^* \rmh)_{s}\big(t, x^{a}, p, p_{b}\big) = \big(s t, x^{a}, p, s p_{b}\big).\end{gather*}
One again this class of induced local coordinates is not necessarily the most convenient. We can make the change $p_a \mapsto \textbf{p}_a = t^{-1} p_a$. The admissible changes of coordinates are then of the form
\begin{alignat*}{3}
& x^{a'} = x^{a'}(x), \qquad && t' = \phi(x)t,& \\
& p' = \phi^{-1}(x) p ,\qquad && \textbf{p}_{b'} = \left(\frac{\partial x^a}{\partial x^{b'}} \right) \left( \textbf{p}_a \frac{\partial}{\partial p} + \frac{\partial }{\partial x^a} \right)p'.&
\end{alignat*}
In this convenient class of local coordinates the base of $\sT^{*}P$ comes with local coordinates $(x^a, p , \textbf{p}_b)$, via the coordinate transformations above we see that $\big(\sT^{*}P \big)_0 = \sJ^1P^+$.

More intrinsically, one can observe that associating with a section $\zs$ of $P^+$ the dif\/feren\-tial~$\xd\zi_\zs$ induces a map which assigns to the f\/irst jet $j^1(\zs)(x)$ of $\zs$ at $x\in P/\Rt$ the class $[\xd\zi_\zs(p)]$ in $\sT^*P/\Rt$, where $p$ is any point of~$P$ projecting to~$x$. This map identif\/ies $J^1P^+$ with $\sT^*P/\Rt$. In local coordinates,
\begin{gather*}\zs(x)=f(x)\ \mapsto \xd\zi_\zs(x,t)=\xd(f(x)t)(x,t)=f(x)\xd t+t\frac{\pa f}{\pa x^a}(x)\xd x^a.
\end{gather*}
The form $f(x)\xd t+t\frac{\pa f}{\pa x^a}(x)\xd x^a$ is $\Rt$-invariant and represents the f\/irst jest $(f(x),\frac{\pa f}{\pa x})$ of~$\zs$ at~$x$.

Note that, since the lifted actions are linear, we have actually a whole series of lifted actions, since the multiplication of $(\sT \rmh)_s$ and $(\sT^*\rmh)_s$ by $s^k$, $k\in\mathbb{Z}$, gives a new principal action. The above one has the advantage that, for $P$ coming from a vector bundle, $P=E^\ti$, it can be extended to the lift of the corresponding action of the monoid $(\R,\cdot)$ of multiplicative reals, which in turn is the most ef\/f\/icient way to obtain the double vector bundle structures on~$\sT E$ and $\sT^* E$~\cite{Grabowski:2006}.

Let now $\zL$ be a Poisson structure on a principal $\mathbb{R}^{\times}$-bundle $(P,\rmh)$, and let~$\zL^\#$ be the corresponding vector bundle morphism
\begin{gather}\label{Ps}
\zL^\#\colon \ \sT^*P\to\sT P.
\end{gather}
The following is straightforward.
\begin{Theorem}\label{liftingactions}
The map \eqref{Ps} intertwines the phase and the tangent lifts of the~$\R^\ti$ action if and only if~$\zL$ is homogeneous of degree~$-1$, i.e.,
$(\rmh_s)_*\zL=s^{-1}\zL$. In this case, we have the reduced map
\begin{gather}\label{reducedmap}\zL^\#_0\colon \ \sT^*P/\Rt\simeq\sJ^1P^+ \lra\sT P/\Rt\simeq\DO^1\big(P^+,P^+\big).
\end{gather}
\end{Theorem}
\begin{Remark}
The above reductions of $\sT P$ and $\sT^*P$ can be extended to a reduction of the whole Courant algebroid structure on the Pontryagin bundle $\sT P\oplus_P\sT^*P$ to a structure of a~\emph{contact Courant algebroid} in the sense of \cite{Grabowski2013} (or \emph{Courant--Jacobi algebroid} in the sense of \cite{Grabowski:2003}) on $\DO^1(P^+,P^+)\oplus_{P_0}\sJ^1P^+$. As we have a canonical pairing
\begin{gather*}\DO^1\big(P^+,P^+\big)\oplus_{P_0}\sJ^1P^+\lra P^+,
\end{gather*}
analogues of Dirac structures can also be naturally def\/ined (see~\cite{Vitagliano:2015}).
\end{Remark}

\section{Kirillov brackets and Kirillov manifolds}\label{sec:Kirillov}
\subsection{Kirillov manifolds} A line bundle equipped with a local Lie bracket on its sections, $(L, [\cdot,\cdot ]_{L})$, as introduced by Kirillov \cite{Kirillov:1976}, is known in the literature also as a \emph{Jacobi bundle} following~\cite{Marle:1991}. Locally this bracket is given by
\begin{gather}\label{jb}[ f,g]_L(x)=\zL^{ab}(x)\frac{\pa f}{\pa x^a}(x)\frac{\pa g}{\pa x^b}(x)+
\zL^a(x)\left(f(x)\frac{\pa g}{\pa x^a}(x)-\frac{\pa f}{\pa x^a}(x)g(x)\right).
\end{gather}
One can identify smooth sections of a line bundle $L$ with smooth homogeneous functions of degree one on $L^*$, and further also with homogeneous functions of degree one on the principal $\R^\ti$-bundle $L^{\*\ti}:=(L^*)^\ti$, i.e., functions $f\colon L^\ti\to\R$ such that $f(\rmh_s(v)):=f(s. v)=s f(v)$. We denote this identif\/ication via $ u \rightsquigarrow \iota_{u}$, where $u \in \Sec(L)$.

Having a Kirillov bracket $[\cdot,\cdot]_L$ on sections of $L$, we can try to def\/ine a Poisson bracket $\{\cdot,\cdot\}_\zL$, associated with a linear Poisson structure $\zL$ on $L^*$, using the identity
\begin{gather}\label{e1}
\iota_{[u, v]_{L}} = \{\iota_{u}, \iota_{v}\}_{\Lambda}.
\end{gather}
However, unlike the case of a Lie algebroid, this bracket is generally singular at points of the zero-section. Instead, one has to def\/ine a Poisson tensor on $L^{*\ti}$. Indeed, in dual coordinates $(x^a,t)$ on $L^*$,
\begin{gather}\label{eqn:possion structure}
\zL(x,t)=\frac{1}{2t}\zL^{ab}(x)\pa_{x^a}\we\pa_{x^b}+\zL^a(x)\pa_t\we\pa_{x^a}.
\end{gather}
This identif\/ication allows for a very useful characterisation of Kirillov brackets (cf.~\cite{Grabowski2013, Marle:1991}) in terms of \emph{Kirillov manifolds} (\emph{Kirillov structures}).

\begin{Definition}[\cite{Grabowski2013}] A \emph{principal Poisson $\R^\ti$-bundle}, shortly \emph{Kirillov manifold}, is a principal $\R^\ti$-bundle $(P,\rmh)$ equipped with a Poisson structure $\zL$ of degree~$-1$. A \emph{morphism of Kirillov manifolds} $\phi\colon P\rightarrow P'$ is a Poisson morphism that intertwines the respective $\R^\ti$-actions.
\end{Definition}

A little more explicitly, a Kirillov manifold is a triple $(P, \rmh, \zL)$, where $(P, \rmh)$ is a principle $\R^\times$-bundle, and $\zL$ is a Poisson structure on $P$ such that $(\rmh_s)_* \zL = s^{-1}\zL$, for all $s \in \R^\times$. In particular, $\R^n\ti\R^\ti$ with coordinates $(x^a,t)$ and trivial $\Rt$-principal bundle structure, equipped with a Poisson tensor of the form (\ref{eqn:possion structure}), is a basic example of a Kirillov manifold.\par

Evidently, Kirillov manifolds form a category under the standard composition of smooth maps. We summarise all the above observations as:

\begin{Theorem}\label{thm:jacobi} There is a one-to-one correspondence between Kirillov brackets $[\cdot,\cdot ]_{L}$ on a line bundle $L \rightarrow M$ and Kirillov manifold structures on the principal $\mathbb{R}^{\times}$-bundle $P = L^{*\ti}$ given by~\eqref{e1}.
\end{Theorem}
\begin{Remark}\label{redu} The local representation (\ref{eqn:possion structure}) identif\/ies (locally) Kirillov manifolds with \emph{Jacobi manifolds} in the sense of Lichnerowicz~\cite{Lichnerowicz:1978}, i.e., manifolds equipped with a bivector f\/ield $\widetilde{\zL} := \frac{1}{2}\zL^{ab}(x)\pa_{x^a}\we\pa_{x^b}$ and a vector f\/ield $\widetilde{R} : = \zL^a(x)\pa_{x^a}$ such that~(\ref{jb}) is a Lie bracket. In terms of the Schouten bracket $[\![\cdot,\cdot]\!]$ this is equivalent to the system of identities
\begin{gather*}[\![\wt{R},\wt{\zL}]\!]=0,\qquad [\![\wt{\zL},\wt{\zL}]\!]=2\wt{R}\we\wt{\zL}.
\end{gather*}

In other words, for a trivial principal $\R^\ti$-bundle $P=M\ti\R^\ti$, the reduced tangent and cotangent bundles can be identif\/ied as
\begin{gather*}
\sT(M\ti\R^\ti)/\R^\ti\simeq \sT M\ti\R,\qquad \sT^*(M\ti\R^\ti)/\R^\ti\simeq \sT^* M\ti\R,
\end{gather*}
and the reduced map (\ref{reducedmap}) is nothing but the vector bundle morphism
\begin{gather*}\zL^\#_0\colon \ \sT^*M\ti\R\ra\sT M\ti\R,\end{gather*}
induced by the Jacobi structure.
\end{Remark}

\subsection{Coisotropic submanifolds of Kirillov manifolds}
Recall that a submanifold $S$ of a Poisson manifold $(P,\zL)$ is called \emph{coisotropic} if the ideal $I_S$ of functions vanishing on $S$ is closed under the Poisson bracket. In the context of a Kirillov manifold $(P, \rmh, \Lambda)$, of particular interest are coisotropic submanifolds which are simultaneously $\Rt$-subbundles; we will call them \emph{coisotropic Kirillov submanifolds} or simply \emph{coisotropic subbundles}. The natural inclusion $S \hookrightarrow P$ implies $S_{0} \hookrightarrow P_0$, where $S_{0} = S \slash \mathbb{R}^{\times}$ is the reduced manifold. Of course, in the traditional language, $S_0$ is called a coisotropic submanifold of the corresponding Jacobi structure (cf.~\cite{Le:2014}).

Suppose now that $S$ is a coisotropic subbundle of the Kirillov manifold $P$. It is clear that in local coordinates $(t, x^{\alpha}, y^{i})$ on $P$ adapted to $S$, so $(t, x^{\alpha})$ form a coordinate system on $S$, the Poisson structure on $P$ encoding the Kirillov manifold structure is of the form
\begin{gather}
\Lambda = \frac{1}{2t} \Lambda^{\alpha \beta}(x,y)\frac{\partial}{\partial x^{\beta}} \wedge\frac{\partial}{\partial x^{\alpha}} + \frac{1}{t}\Lambda^{\alpha i}(x,y)\frac{\partial}{\partial y^{i}} \wedge \frac{\partial}{\partial x^{\alpha}}\nonumber\\
\hphantom{\Lambda =}{} + \frac{1}{2t} \Lambda^{ij}(x,y)\frac{\partial}{\partial y^{j}} \wedge\frac{\partial}{\partial y^{i}} + \Lambda^{\alpha}(x,y)\frac{\partial}{\partial x^{\alpha}} \wedge\frac{\partial}{\partial t} + \Lambda^{i}(x,y)\frac{\partial}{\partial y^{i}} \wedge\frac{\partial}{\partial t}, \label{cps}
\end{gather}
where we require $\Lambda^{ij}=0$ and $\Lambda^{i}=0$ on $S$.

Hamiltonian vector f\/ields $X_{f} := \{ f, \cdot\}_{\Lambda}$ where $f \in I_{S}$ are tangent to $S$ and so form an integrable distribution. The corresponding foliation is known as the \emph{characteristic foliation} of~$S$. The space of leaves, if smooth, will inherit a Kirillov manifold structure. The reduction of the Poisson structure~(\ref{cps}) is obvious. To see that we can canonically reduce the homogeneous structure, consider the weight vector f\/ield $\nabla_P$ generating $h$ and the Lie bracket $[\nabla_P,X_f]$, where $f\in I_S$. We have
\begin{gather*} \lan[\nabla_P,X_f],\xd g\ran= \Ll_{\nabla_P}\lan {X_f},{\xd g}\ran-\lan X_f,\Ll_{\nabla_P}\xd g\ran\\
\hphantom{\lan[\nabla_P,X_f],\xd g\ran}{} = \Ll_{\nabla_P}\left(\zL(\xd f,\xd g)\right)-\zL(\xd f,\Ll_{\nabla_P}\xd g)\\
\hphantom{\lan[\nabla_P,X_f],\xd g\ran}{} = (\Ll_{\nabla_P}\zL)(\xd f,\xd g) +\zL(\Ll_{\nabla_P}\xd f,\xd g)\\
\hphantom{\lan[\nabla_P,X_f],\xd g\ran}{} = -\zL(\xd f,\xd g)+\zL(\Ll_{\nabla_P}\xd f,\xd g)=X_{(\nabla_P(f)-f)}(g).
\end{gather*}
Here we used the fact that $\zL$ is homogeneous. Now, as $\nabla_P$ is tangent to $S$, we have \smash{$\nabla_P(f){-}f\!\in\! I_S$}, so that~$\nabla_P$ lies in the normalizer of the Lie algebra of vector f\/ields tangent to leaves of the characteristic distribution, hence~$\nabla_P$ (as well as~$\rmh$) preserves the foliation and thus induces a~homogeneous structure on the manifold of leaves.

Studying reductions, deformations, etc., of coisotropic subbundles is an interesting but extensive task, which we postpone for future study.

\subsection{Contact structures}
Proceeding to contact structures, f\/irst note that a nowhere-vanishing one-form $\za$ spans a trivial one dimensional vector subbundle $[\alpha]$ of $\sT^{*}M$. Associated with $\alpha$ is a canonical embedding $I_{\alpha}\colon \mathbb{R} \times M \rightarrow \sT^{*}M$ which induces an isomorphism of $\mathbb{R} \times M$ with $[\alpha]$. In natural local coordinates the canonical embedding is given by
\begin{gather}\label{lc}
I_{\alpha}^{*}\big(x^{a}, p_{b}\big) = \big(x^{a}, t \alpha_{b}(x)\big),
\end{gather}
where $t$ is the (global) coordinate on $\mathbb{R}$ and locally we have $\alpha = \alpha_{a}(x) \rmd x^{a}$.
\begin{Proposition}
The nowhere-vanishing one-form $\alpha$ is a contact form if and only if the trivial principal bundle $\mathbb{R}^{\times}\times M$ is, via $I_{\alpha}$, a symplectic submanifold $[\alpha]^{\times} \subset \sT^{*}M$.
\end{Proposition}

The above propositions is essentially a well-known rewording of the standard notion of the `symplectisation' of a contact form. In particular, it is easy to see that in Darboux coordinates
\begin{gather*}
I_{\alpha}^{*}\big(\rmd p_{a} \wedge \rmd x^{a}\big) := \omega = \rmd t\wedge \alpha , 
\end{gather*}
\noindent which gives the symplectisation of $\alpha$ remembering that $t \neq 0$. Moreover, note that the contact form can be recovered from
\begin{gather*}
i_{\nabla}\omega = t \alpha,
\end{gather*}
\noindent where $\nabla$ is the Euler vector f\/ield on the principal $\mathbb{R}^{\times}$-bundle $[\alpha]^{\times}$, i.e., the fundamental vector f\/ield of the $\Rt$-action. All this implies the following.
\begin{Proposition}[\cite{Grabowski2013}] A line subbundle $C$ of $\sT^{*}M$ is locally generated by contact one-forms if and only if $C^{\times}$ is a symplectic submanifold of $\sT^{*}M$.
\end{Proposition}
\begin{Definition}
A principal $\mathbb{R}^{\times}$-bundle $(P, \rmh)$ equipped with a 1-homogeneous symplectic form~$\omega$, i.e., a symplectic form such that $(\rmh_{t})^{*}\omega = t \omega$ ($t \neq 0$), will be referred to as a \emph{contact structure}. In other words, a contact structure is a Kirillov manifold whose Poisson structure is invertible (symplectic).
\end{Definition}
Let $\nabla$ be the Euler vector f\/ield on $P$, $\nabla\colon P\to\sT P$. It is easy to see that the composition $\zh=\zw^\flat\circ\nabla\colon P\to\sT^*P$ is a one-form on $P$ which takes values in basic covectors, $\zh(y)=\zp^*(\Psi(y))\in\sT^*_yP$, $\Psi(y)\in\sT^*_{\zp(y)}P_0$, so can be viewed as a map $\Psi\colon P\to\sT^*P_0$ which locally has the form~(\ref{lc}). Consequently, the range $C^\ti(P)=\{\Psi(y)~|~y\in P\}$ of $\Psi$ is a a symplectic submanifold in $\sT^*P_0$. Thus we get the following.
\begin{Theorem}[\cite{Grabowski2013}]\label{thm:contact} Any contact structure $(P,\zw,\rmh)$, where $P$ is an $\mathbb{R}^{\times}$-bundle over~$P_0$, can be canonically symplectically embedded into $\sT^{*}P_0$ as a symplectic principal $\mathbb{R}^{\times}$-bundle of the form~$C^{\times}$ for a~line subbundle $C\subset\sT^*P_0$.
\end{Theorem}

\begin{Remark}\label{rem:contact distribution} Commonly, a contact structure on a manifold $M$ is understood as a maximally non-integrable hyperplane distribution $\mathcal{D} \subset \sT M$, locally given as the polar (annihilator) of a line bundle $C \subset \sT^{*} M$ generated by contact one-forms, $\mathcal{D}=C^0$. We will refer to such hyperplane distributions $\mathcal{D} \subset \sT M$ as \emph{contact distributions} to avoid confusion. In our language, contact structures are homogeneous symplectic structures on a principal $\R^\ti$-bundle $P$, while in the classical language they are certain hyperplane distributions on the reduced manifold $M=P/\R^\ti$.
\end{Remark}

\begin{Example} The canonical symplectic structure on the cotangent bundle $\sT^*M$ is linear, thus homogeneous on $(\sT^*M)^\ti$. The symplectic homogeneous manifold $P=(\sT^*M)^\ti$ represents therefore a contact structure. In the traditional language it is a canonical contact structure on the reduced manifold $(\sT^*M)^\ti/\Rt$, i.e., on the projectivisation $\cP(\sT^*M)$ of the cotangent bundle.
\end{Example}

\begin{Example} Consider a principal $\R^\ti$-bundle $(P,\rmh)$. It is easy to see that the canonical symplectic form on the cotangent bundle $\sT^* P$ is homogeneous with respect to the lifted action~$\sT^* \rmh$, so~$\sT^* P$ represents canonically a contact structure. If we write $P=L^{*\ti}$, then in the traditional language this is exactly the canonical contact structure $C$ on the reduced mani\-fold~$\sT^*P/\R^\ti$ which is the f\/irst jet bundle $J^1P^+=J^1L$. When $L$ is the trivial bundle, i.e., $L = \R \times M$, the canonical contact structure $C$ is the trivial line subbundle of $\sT^{*}(\mathbb{R} \times \sT^{*}M)$ generated by the contact form $\alpha = \rmd z - p_{a} \rmd x^{a}$. Thus we have $P=C^{\times} = \mathbb{R}^{\times} \times \mathbb{R} \times \sT^{*}M$, which we equip with local coordinates $(t, z, x^{a}, p_{b})$ and thus the symplectic structure on $P$ is
\begin{gather*}
\omega = \rmd t \wedge \rmd z - p_{a}\rmd t \wedge \rmd x^{a} - t \rmd p_{a} \wedge \rmd x^{a}.
\end{gather*}
\end{Example}

\begin{Remark}Contact structures on non-negatively graded manifolds further equipped with homological contact vector f\/ields were studied by Mehta~\cite{Mehta:2013} using a more traditional language than put forward here. In particular, for the degree~$1$ case he showed that there is a~one-to-one correspondence between such structures (with a global contact form) and Jacobi manifolds. The line bundle approach to the concept of a~\emph{generalised contact bundle} can be found in the work of Vitagliano and Wade~\cite{Vitagliano:2016}. Furthermore, the $\Rt$-principal bundle approach can also be applied to the notion of a \emph{contact structure on a Lie algebroid} following Ida and Popescu \cite[Remark~4.2]{Ida:2016}.
\end{Remark}

\section{Principal bundle Lie groupoids and algebroids}\label{sec:Ggroupoids}

\subsection{Morphisms of Lie groupoids and Lie algebroids}
Our general reference to the theory of Lie groupoids and Lie algebroids will be Mackenzie's book~\cite{Mackenzie2005}.

Let $\mathcal{G} \rightrightarrows M$ be an arbitrary Lie groupoid with \emph{source map} $s\colon \mathcal{G} \rightarrow M$ and \emph{target map} $t\colon \mathcal{G} \rightarrow M$. There is also the inclusion map $\iota_M \colon M \rightarrow \mathcal{G}$, $\zi_M(x)=\Id_x$, and a~\emph{partial multiplication} $(g,h) \mapsto gh$ which is def\/ined on $\mathcal{G}^{(2)}=\{(g,h)\in\mathcal{G}\ti\mathcal{G}\colon s(g) = t(h)\}$. Moreover, the mani\-fold~$\mathcal{G}$ is foliated by $s$-f\/ibres $\mathcal{G}_{x}= \{ \left.g \in \mathcal{G}\right| s(g) =x\}$, where $x \in M$. As by def\/inition the source and target maps are submersions, the $s$-f\/ibres are themselves smooth manifolds. Geometric objects associated with this foliation will be given the superscript~$s$. In particular, the distribution tangent to the leaves of the foliation will be denoted by $\sT^{s} \mathcal{G}$. To ensure no misunderstanding with the notion of a Lie groupoid morphism we recall the def\/inition we will be using.

\begin{Definition}\label{def:Lie groupoid Morphism}
Let $\mathcal{G}_{i} \rightrightarrows M_{i}$ $(i=1,2)$ be a pair of Lie groupoids. Then a \emph{Lie groupoid morphisms} is a pair of maps $(\Phi, \phi)$ such that the following diagram is commutative
\begin{gather*}
\begin{xy}
(0,20)*+{\mathcal{G}_{1}}="a"; (20,20)*+{\mathcal{G}_{2}}="b";%
(0,0)*+{M_{1}}="c"; (20,0)*+{M_{2}}="d";%
{\ar "a";"b"}?*!/_2mm/{\Phi};
{\ar@<1.ex>"a";"c"} ;
{\ar@<-1.ex> "a";"c"} ?*!/^3mm/{s_{1}} ?*!/_6mm/{t_{1}};
{\ar@<1.ex>"b";"d"};%
{\ar@<-1.ex> "b";"d"}?*!/^3mm/{s_{2}} ?*!/_6mm/{t_{2}}; %
{\ar "c";"d"}?*!/^3mm/{\phi};
\end{xy}
\end{gather*}
in the sense that $s_{2}\circ \Phi = \phi \circ s_{1}$, and $t_{2}\circ \Phi = \phi \circ t_{1} $
subject to the further condition that $\Phi$ respects the (partial) multiplication; if $g,h \in \mathcal{G}_{1}$ are composable, then $ \Phi(gh) = \Phi(g)\Phi(h)$. It then follows that for $x \in M_{1}$ we have $\Phi(\Id_{x}) = \Id_{\phi(x)}$ and $\Phi(g^{-1}) = \Phi(g)^{-1}$.
\end{Definition}

Consider a Lie groupoid $\mathcal{G} \rightrightarrows M$. A \emph{Lie subgroupoid} of $\mathcal{G}$ is a Lie groupoid $\mathcal{H}\rightrightarrows M'$ together with injective immersions $\zi\colon \mathcal{H} \to \mathcal{G}$ and $\zi_0\colon M'\to M$ such that $(\zi,\zi_0)$ is a morphism of Lie groupoids. In particular, $\zi(\mathcal{H})$ it is closed under multiplication (when def\/ined) and inversion. If $\zi_0(M')=M$, then the Lie subgroupoid $\mathcal{H}\rightrightarrows M$ is said to be a \emph{wide subgroupoid}.

The Cartesian product $\mathcal{G}_{1}\ti\mathcal{G}_{2}$ of two Lie groupoids is canonically a Lie groupoid, and it follows immediately from the above def\/inition that $\Phi$ is a Lie groupoid morphism if and only if its graph is a Lie subgroupoid in $\mathcal{G}\ti\mathcal{H}$.

A similar fact holds true for Lie algebroids, but as we have many alternative def\/initions of a Lie algebroid, there are many alternative def\/initions of a Lie algebroid morphism (see, e.g., \cite[Theorem~14]{Grabowski2012a} in a little more general setting). Dealing with homogeneous Poisson structures in this paper, we will mainly understand a Lie algebroid on a vector bundle $E$ as a linear Poisson structure on~$E^*$. Then, as is commonly known (see, e.g., Mackenzie \cite[p.~400]{Mackenzie2005}), Lie subalgebroids of $E$ correspond to coisotropic subbundles in the Poisson manifold $E^*$ by passing to the polar in the dual bundle.

It is also well known that via a dif\/ferentiation procedure one can construct the \emph{Lie functor}
\begin{gather*}
\catname{Grpd} \stackrel{\Lie}{\xrightarrow{\hspace*{30pt}}} \catname{Algd},
\end{gather*}
that sends a Lie groupoid to its Lie algebroid, and sends morphisms of Lie groupoids to morphisms of the corresponding Lie algebroids. However, as is also well known, we do not have an equivalence of categories as not all Lie algebroids arise as the inf\/initesimal versions of Lie groupoids. There is no direct generalisation of Lie III, apart from the local case. The obstruction to the integrability of Lie algebroids, the so called \emph{monodromy groups}, were f\/irst uncovered by Crainic and Ferandes \cite{Crainic:2003}. To set some notation and nomenclature, given a Lie groupoid $\mathcal{G}$, we say that $\cG$ \emph{integrates} $\Lie(\mathcal{G}) = \rmA(\mathcal{G})$. Moreover, if $\Phi\colon \mathcal{G} \rightarrow \mathcal{H}$ is a morphism of Lie groupoids, then we will write $\Phi' = \Lie(\Phi) \colon \rmA(\mathcal{G}) \rightarrow \rmA(\mathcal{H})$ for the corresponding Lie algebroid morphism, which actually comes from the dif\/ferential $\sT\Phi\colon \sT\mathcal{G} \rightarrow \sT\mathcal{H}$ restricted to the $s$-f\/ibres.

Let us just recall Lie~II theorem as we will need it later on.
\begin{Theorem}[Lie II]\label{thm:integration morphisms} Let $\mathcal{G} \rightrightarrows M$ and $\mathcal{H}\rightrightarrows N$ be Lie groupoids. Suppose that $\mathcal{G}$ is source simply-connected and that $\phi\colon \rmA(\mathcal{G}) \rightarrow \rmA(\mathcal{H})$ is a Lie algebroid morphism between the associated Lie algebroids. Then, $\phi$ integrates to a unique Lie groupoid morphisms $\Phi\colon \mathcal{G} \rightarrow \mathcal{H}$ such that $\Phi' = \phi$.
\end{Theorem}

This generalisation of Lie II to the groupoid case was f\/irst proved by Mackenzie and Xu~\cite{Mackenzie:2000}. A simplif\/ied proof was obtained shortly after by Moerdijk and Mr\v{c}un~\cite{Moerdijk:2002}. Note that the assumption that the Lie groupoid $\mathcal{G}$ is source simply-connected is essential.

\subsection{Compatible group actions on Lie groupoids and algebroids}
 In our study of Jacobi and contact groupoids we will encounter Lie groupoids that have a~compa\-tible action of $\mathbb{R}^{\times}$ upon them; compatibility to be def\/ined shortly. However, as the basic theory of compatible group actions on Lie groupoids is independent of the actual Lie group, we discuss the general setting here focusing on what we will need later in this paper.

\begin{Definition}\label{d1} An action $\rmh\colon G\ti\cG\to\cG$ of a Lie group $G$ on a Lie groupoid $\cG\rightrightarrows M$ is said to be \emph{compatible with the groupoid structure} if $\rmh_g\colon \cG\to\cG$ are groupoid isomorphisms for all $g\in G$. A principal $G$-bundle $\zp\colon \cG\to\cG_0$ is a \emph{principal bundle $G$-groupoid} (\emph{$G$-groupoid in short}) if the principal action of $G$ on $\cG$ is compatible with the groupoid structure. Similarly, a $G$-action on a Lie algebroid $A$ is \emph{compatible} if the group acts by Lie algebroid isomorphisms, and we get a~\emph{$G$-algebroid} if a principal $G$-action is compatible with the Lie algebroid structure.
\end{Definition}

\begin{Remark}The notion of a group object in the category of groupoids is quite an old notion, going back at least to the mid 1960s with the unpublished works of Verdier and Duskin. The earliest published work on such objects that we are aware of is that of Brown and Spencer~\cite{Brown:1976}.
\end{Remark}

\begin{Remark}The reader should also be reminded of Mackenzie's notion of a~\emph{PBG-groupoid} \cite{Mackenzie:1987,Mackenzie:1988}, which is close to our notion of a G-groupoid, although Mackenzie, being interested in extensions of principal bundles, starts with a principal $G$-structure on the manifold $M$ of units extended accordingly to a Lie groupoid~$\cG$. The other dif\/ference is that what we call a Lie groupoid is a~\emph{differentiable groupoid} in the sense of Mackenzie, and his \emph{Lie groupoids} in~\cite{Mackenzie:1987,Mackenzie:1988} (or \emph{locally trivial Lie groupoids} in~\cite{Mackenzie2005}) form much smaller class and are understood as particular \emph{transitive Lie groupoids}. However, the following observations are independent of these details and so are probably already known to Mackenzie.
\begin{enumerate}\itemsep=0pt
\item The action of $G$ on $\mathcal{G}$ commutes with the source and target maps, thus projects onto a~$G$-action on the manifold $M$. Moreover, $M$ as an immersed submanifold of $\cG$ is invariant with respect to the $G$-action, and the projected and restricted actions coincide.
\item As the action of $G$ on $\mathcal{G}$ is principal, it is also principal on the immersed submani\-fold~$M$, so $M$ inherits a structure of a principal $G$-bundle. It is important to note that $M$ is G-invariant. In particular, the quotient manifold $M_0=M/G$ exists.
 \item The reduced manifold $\cG/G=\mathcal{G}_{0}$ is a Lie groupoid $\cG/G=\mathcal{G}_{0} \rightrightarrows M/G=M_{0}$, with the set of units $M_0$, def\/ined by the following structure:
\begin{gather*}
\begin{split}
& \begin{xy}
(0,20)*+{\mathcal{G}}="a"; (20,20)*+{\mathcal{G}_{0}}="b";%
(0,0)*+{M}="c"; (20,0)*+{M_{0}}="d";%
{\ar "a";"b"}?*!/_2mm/{\pi};
{\ar@<1.ex>"a";"c"} ;
{\ar@<-1.ex> "a";"c"} ?*!/^3mm/{s} ?*!/_6mm/{t};
{\ar@<1.ex>"b";"d"};%
{\ar@<-1.ex> "b";"d"}?*!/^3mm/{{\zs}} ?*!/_6mm/{{\zt}}; %
{\ar "c";"d"}?*!/^3mm/{\rmp};
\end{xy}\end{split}
\qquad
\begin{array}{l}
 {\zs} \circ \pi = \rmp \circ s, \\
 {\zt} \circ \pi = \rmp \circ t,\\
 \Id_{\rmp(x)}=\zp(\Id_x)\quad\text{for all}\quad x\in M,\\
 \pi(y)^{-1}=\pi\big(y^{-1}\big)\quad\text{for all}\quad y\in\cG,\\
 \pi(y y') = \pi(y) \pi(y')\quad\text{for all}\quad(y,y') \in \mathcal{G}^{(2)},
\end{array}
\end{gather*}
where $\pi\colon \cG\to\cG_0$ is the canonical projection. The source map~$\zs$ is clearly a submersion. Indeed, as~$\rmp$ and $s$ are submersions, their composition $\rmp\circ s=\sigma\circ\pi$ is a submersion and therefore~$\sigma$ is also a submersion.
\end{enumerate}

In fact, the above constructions imply, tautologically, that $(\pi, \rmp) \colon \mathcal{G} \rightrightarrows M \rightarrow \mathcal{G}_{0} \rightrightarrows M_{0}$ is a~morphism of Lie groupoids with the above structures.

The concept of a $G$-groupoid is essentially of double nature: a $G$-groupoid is a principal $G$-bundle object in the category of Lie groupoids. From the point of view of Jacobi and contact geometry, the most important will be of course $\R^\ti$-groupoids.
\end{Remark}

\begin{Remark} It is well known that a $G$-action on a set $X$ is equivalent to a groupoid morphism of $G$ into the pair groupoid $X\ti X$. It can be shown that if $X=\cG$ is a (Lie) groupoid and the action is by automorphism, then the morphism of $G$ into $\cG\ti\cG$ is simultaneously a morphism
with respect to the other, namely Cartesian product groupoid structure on $\cG\ti\cG$ (with $M\ti M$ as the set of units, which is simultaneously the pair groupoid over $M$).
A $G$-groupoid can be therefore also def\/ined as a (double) groupoid morphism of $G$ (viewed as a double groupoid) into the double groupoid (in the sense of Ehresmann) $\cG\ti\cG$, with the diagram
\begin{gather*}
\begin{xy}
(0,25)*+{\mathcal{G}\ti\cG}="a"; (25,25)*+{\mathcal{G}}="b";%
(0,0)*+{M\ti M}="c"; (25,0)*+{M.}="d";%
{\ar@<1.ex>"a";"b"} ;
{\ar@<-1.ex> "a";"b"}; 
{\ar@<1.ex>"a";"c"} ;
{\ar@<-1.ex> "a";"c"} ?*!/^6mm/{s\ti s} ?*!/_8mm/{t\ti t};
{\ar@<1.ex>"b";"d"};%
{\ar@<-1.ex>"b";"d"}?*!/^3mm/{s} ?*!/_6mm/{t}; %
{\ar@<1.ex>"c";"d"};
{\ar@<-1.ex>"c";"d"} ;
\end{xy}
\end{gather*}
 We are unable to fully investigate the corresponding theory here: replacing the group $G$ with a~groupoid leads to groupoid morphisms in the sense of Zakrzewski \cite{Stachura:2000,Zakrzewski:1990a}, which are nowadays also called \emph{groupoid comorphisms}.
\end{Remark}

It is easy to see that a compatible principal $G$-structure on a Lie groupoid $\cG$ induces
canonically a compatible principal $G$-structure on the Lie algebroid $\Lie(\cG)$. Indeed, if $\rmh\colon G\ti\cG\to\cG$ is such a structure, then \emph{via} the f\/irst Lie theorem,
\begin{gather}\label{lie}
\rmh'_g=\Lie(\rmh_g)
\end{gather}
def\/ines a free $G$-action on $\Lie(\cG)$ by automorphisms. This action is also proper, as the Lie functor is a restriction of the tangent functor; the tangent lift of a proper group action is a~proper group action. Actually we have the following theorem on integrability of $G$-algebroids.

 \begin{Theorem}\label{thm:lie theory01}
There is a one-to-one correspondence between compatible $G$-structures on $\Lie(\mathcal{G})$ and on $\mathcal{G}$ satisfying \eqref{lie}, for $\mathcal{G}$ source simply connected.
\end{Theorem}

\begin{proof}
It remains to prove that if $g\mapsto \rmh'_g$ gives rise to a compatible principal $G$-action on $\Lie(G)$, then $g\mapsto \rmh_g$ is also principal and compatible. To check that the $G$-action on $\mathcal{G}$ is principal, one can use the general fact about Poisson actions of Poisson--Lie groups proven in \cite[proof of Proposition~3.1]{Fernandes:2009}. Here, for the convenience of the reader we present a direct proof.

Via Lie II we know that the latter is a free $G$-group action as Lie groupoid automorphisms. This fact that this action is smooth also follows from Lie II. In particular, applying Lie II to $G \times \Lie(\cG) \rightarrow \Lie(\cG)$ (considering $G$ as the base of a rank zero Lie algebroid) yields a smooth map $G \times \cG \rightarrow \cG$ (now thinking of $G$ as the unit groupoid over itself). The fact that $\rmh_{g_1 g_2} = \rmh_{g_1} \circ \rmh_{g_2}$ follows from the uniqueness of integration applied to each $\rmh_g$ separately. It only remains to show that this group action is proper.

Of course, $M$ is an invariant submanifold of this action and the `integrated' action coincides with the original action on $M\subset \Lie(\cG)$, thus is proper. Moreover, the integrated action on $\mathcal{G}$ projects \emph{via} the source map $s$ to the action on $M$ which implies that the integrated action is proper. Indeed, having two compact sets $K_i$, $i=1,2$, in $\cG$, we have that \smash{$\{ g\in G\, |\, \rmh_g(K_1)\cap K_2\ne\varnothing\}$} is a closed subset of the compact set $\{ g\in G\, |\, \rmh_g(s(K_1))\cap s(K_2)\ne\varnothing\}$, thus compact.
\end{proof}

\begin{Remark}
The above theorem can be derived from the main results of Stefanini \cite{Stefanini:2007,Stefanini:2008} describing, roughly speaking, integrability conditions for $G$-algebroids with $G$ being a Lie groupoid. However, the Lie group case is substantially simpler, so we decided to present the direct proof.
\end{Remark}

\subsection[Structure of $G$-groupoids]{Structure of $\boldsymbol{G}$-groupoids}
Let now $\cG$ be a $G$-groupoid with the structure diagram
\begin{gather}\label{pGgr}\begin{split}& \begin{xy}
(0,20)*+{\mathcal{G}}="a"; (20,20)*+{\mathcal{G}_{0}}="b";%
(0,0)*+{M}="c"; (20,0)*+{M_{0}.}="d";%
{\ar "a";"b"}?*!/_2mm/{\pi};
{\ar@<1.ex>"a";"c"} ;
{\ar@<-1.ex> "a";"c"} ?*!/^3mm/{s} ?*!/_6mm/{t};
{\ar@<1.ex>"b";"d"};%
{\ar@<-1.ex> "b";"d"}?*!/^3mm/{{\zs}} ?*!/_6mm/{{\zt}}; %
{\ar "c";"d"}?*!/^3mm/{\rmp};
\end{xy}\end{split}
\end{gather}
\begin{Proposition} The map
\begin{gather}\label{sm0}
\cS\colon \ \cG\to \cG_0\times_{M_0}M:= \big\{(y_0,x) \in \mathcal{G}_{0}\ti M \, | \, \rmp(x)= {\zs}(y_0) \big\},\qquad
\cS(y)= (\zp(y),s(y) ),
\end{gather}
is a dif\/feomorphism. This dif\/feomorphisms identifies $\cG$ as a $G$-bundle over $\cG_0$ with the pull-back of $\rmp\colon M\to M_0$ along $\zs\colon \cG_0\to M_0$ which is $\rmp^!\cG_0=\cG_0\times_{M_0}M$.
\end{Proposition}
\begin{proof}
The map ${\cS}$ is clearly a dif\/feomorphism if the principal $G$-bundle $\cG$ is trivial. Indeed, if $\cG=\cG_0\ti G$, then $M\simeq M_0\ti G$, the manifold $\cG_0\times_{M_0}M$ can be identif\/ied with $\cG_0\ti G$ by $(y_0,g)\mapsto(y_0,({\zs}(y_0),g))$, and with these identif\/ications the map ${\cS}\colon \cG_0\ti G\to \cG_0\ti G$ is the identity. As the bundle $\cG\to\cG_0$ is locally trivial, the map ${\cS}$ is generally a surjective local dif\/feomorphism. It is also globally injective, thus a global dif\/feomorphism. Indeed, ${\cS}(y)={\cS}(y')$ implies that $\zp(y)=\zp(y')$, so $y'=yg$ for some $g\in G$, and therefore $s(y')=s(y)g$. But ${\cS}(y)={\cS}(y')$ implies also $s(y')=s(y)$, so that $g=e$ (the action is free) and $y=y'$. As, clearly, $(y_0,x)g=(y_0,xg)$, the $G$-bundle $\cG$ is the pull-back bundle.
\end{proof}

Using the above identif\/ication, we can transmit the $G$-groupoid structure from $\cG$ onto \smash{$\cG_0\times_{M_0}M$}. The $G$-action is clearly $(y_0,x)g=(y_0,xg)$, the embedding of units is $\zi_M(x)=(\Id_x,x)$, and the source map reads $s(y_0,x)=x$. Knowing the inverse we could def\/ine the target map and the composition. It is easy to see that the inverse of $y=(y_0,x)$ is $y^{-1}=(y_0^{-1},t(y_0,x))$, where $t$ is the target map. However, the groupoid structure on $\cG$ is not the pull-back groupoid structure corresponding to the pull-back of $\cG_0$ along $\rmp\colon M\to M_0$. The argument is dimensional: $\dim(\cG)=\dim(\cG_0)+\dim(G)$, while the pull-back groupoid
\begin{gather*}\{(x,y_0,x')\in M\ti\cG_0\ti M\,|\, \rmp(x)=\zt(y_0),\ \rmp(x')=\zs(y_0)\}\end{gather*}
is of dimension $\dim(\cG_0)+2\dim(G)$.

One can easily check what properties of $t$ ensure that the axioms of a groupoid hold true.
\begin{Theorem}\label{generalsplit} Let $\rmp\colon M\to M_0$ be a principal $G$-bundle with the right $G$-action $M\ti G\ni(x,g)\mapsto xg\in M$, and $\cG_0\rightrightarrows M_0$ be a Lie groupoid with the source and the target map ${\zs}$ and ${\zt}$, respectively.

Then, any $G$-groupoid structure on the manifold $\cG_0\times_{M_0}M$ equipped with the principal $G$-action $(y_0,x)g=(y_0,xg)$, the source map $s(y_0,x)=x$, and such that the projection $(y_0,x)\mapsto y_0$ is a groupoid morphism, is uniquely determined by its target map~$t$. On the other hand, a map $t\colon \cG_0\times_{M_0}M\to M$, $t(y_0,x)=:y_0. x$, can serve as such a target map if and only if it has the following properties $($holding for all $x\in M)$:
\begin{itemize}\itemsep=0pt
\item[$(i)$] $\rmp(y_0.x)=\zt(y_0)$ for all $y_0\in\cG_0$,
\item[$(ii)$] $y_0.(y'_0.x)=(y_0y'_0).x$ for all $(y_0,y'_0)\in\cG_0^2$,
\item[$(iii)$] $\Id_{\rmp(x)}.x=x$,
\item[$(iv)$] $y_0.(xg)=(y_0.x)g$ for all $y_0\in\cG_0$ and all $g\in G$.
\end{itemize}
\end{Theorem}
In simple terms, $(i)$--$(iii)$ mean that $t$ is an action of $\cG_0$ on $\rmp\colon M\to M_0$ (cf.\ \cite[De\-f\/i\-ni\-tion~1.6.1]{Mackenzie2005}), and $(iv)$ means that the action is $G$-equivariant. The $G$-groupoid determined by $t$ as above we will denote $\cG_0\times_{M_0}^t M$ and called \emph{$t$-fixed $G$-groupoid}. Thus, any $G$-groupoid (\ref{pGgr}) is $t$-f\/ixed for some $t(y_0,x)=y_0.x$ satisfying $(i)$--$(iv)$.

There are two particular cases of the above construction which are of great importance. The f\/irst is the case of a trivial principal bundle, $M=M_0\ti G$ which is always a local form of any $G$-groupoid. In this case we can use the identif\/ication $\cG_0\times_{M_0}M\simeq\cG_0\times G$ and replace the map~$t$ satisfying~$(i)$ with a map $\rmb\colon \cG_0\to G$. Indeed, any map on a Lie group commuting with all the right-translations is a left-translation, so can we write $t(y_0,\zs(y_0),g)=(\zt(y_0),\rmb(y_0)g)$. Now, the properties $(i)$--$(iv)$ can be reduced to
\begin{gather*}
\rmb(y_0)\rmb(y'_0)=\rmb(y_0y'_0)
\end{gather*}
for all $(y_0,y'_0)\in\cG_0^2$, i.e., to the fact that $\rmb\colon \cG_0\to G$ is a groupoid morphism. This is of course always the local form of any $G$-groupoid. The corresponding $G$-groupoid structure, denoted with $\cG_0\ti^\rmb G$, is an obvious generalisation of the groupoid extension by the additive $\R$ with a~help of a~multiplicative function considered in the literature (cf.~\cite[Def\/inition~2.3]{Crainic:2007}), and we have shown that this construction is in a~sense universal. Thus we get the following.
\begin{Theorem}\label{trivialsplit} For any $G$-groupoid structure on the trivial $G$-bundle $\cG=\cG_0\ti G$ there is a~Lie groupoid structure on $\cG_0$ with the source and target maps $\zs,\zt\colon \cG_0\to M_0$ and a groupoid morphism $\rmb\colon \cG_0\to G$ such that the source map~$s$, the target map $t$ and the partial multiplication in $\cG$ read
\begin{gather*}s(y_0,g)=(\zs(y_0), g),\qquad t(y_0,g)=(\zt(y_0),\rmb(y_0)g),\qquad (y_0,g_1)(y'_0,g_2)=(y_0y'_0,g_2).
\end{gather*}
\end{Theorem}

Another particular case is that of a bundle of groups, i.e., a groupoid in which the source and the target map coincide (the anchor map $\zr=(s,t)$ is \emph{diagonal}). This means that $\cG_0$ is a~bundle of groups and the map~$t$ is trivial, $t(y_0,x)=x$. Any $G$-groupoid with diagonal anchor splits therefore as the product $\cG=\cG_0\times_{M_0}M$ in which all groupoid operations come from $\cG_0$ and the principal $G$-action from $M$.

This is in particular the case of a $G$-vector bundle, i.e., a vector bundle $\zt\colon P\to M$ on which $G$ acts principally by vector bundle automorphisms, which means in this case that the $G$-action commutes with the natural homogeneity structure ${l} \colon \mathbb{R} \times P \rightarrow P$ that is associated with homotheties of the said vector bundle structure. In other words, $(tv)g=t(vg)$ and we have the diagram
\begin{gather*}
\xymatrix{
P\ar[rr]^{\zp} \ar[d]^{\zt} && P_0\ar[d]^{{\zt_0}} \\
M\ar[rr]^{{\zp_0}} && M_0, }
\end{gather*}
where, as we already know, $\zt$, $\zt_0$ are vector bundles, and $\zp$, $\zp_0$ are principal $G$-bundles. We stress that, even when $G =\R^\ti$, this double structure is not a double vector bundle. In particular, $P_{0}$~is \emph{not} canonically embedded in~$P$, but we have a variant of isomorphism~(\ref{sm0}),
\begin{gather}\label{ident}
(\zt,\zp)\colon \ P\to M\ti_{M_0}P_0=\{(x,y_0)\in M\ti P_0\colon \zp_0(x)=\zt_0(y_0)\}.
\end{gather}
In other words, we get the following generalisation of \cite[Theorem~3.2]{Grabowski2013}.
\begin{Theorem} If $\zt\colon P\to M$ is a $G$-vector bundle, then there is an induced principal $G$-action on $M$ and a splitting $E=M\ti_{M_0}P_0$, where $P_0=P/G$ is a vector bundle over $M_0=M/G$. The $G$-vector bundle structure on $P$ comes directly from this splitting in the obvious manner.
\end{Theorem}

\subsection[Linear $\R^\ti$-bundles]{Linear $\boldsymbol{\R^\ti}$-bundles}
When dealing with Jacobi and contact geometry, principal $\mathbb{R}^{\times}$-bundles $\pi \colon P\rightarrow P_{0} $ that also carry a compatible vector bundle structure $\zt\colon P\to M$ are an essential part of the theory. We will refer to such structures as \emph{linear $\R^\ti$-bundles}.

We are free to employ local homogeneous coordinates of the form $(t, x^{\alpha}, y^{i})$ on $P$, where $(t, x^{\alpha})$ represent coordinates on $M$ and $(x^{\alpha}, y^{i})$ on $P_{0}$, so that the $\R^\ti$-action $\rmh$ reads
\begin{gather*}
\rmh_s\big(t, x^{\alpha}, y^{i}\big)=\big(s t, x^{\alpha}, y^{i}\big)
\end{gather*}
 and $(t, x^{\alpha}, y^{i})\mapsto (t, x^{\alpha})$ is a vector f\/ibration.

We will use the following fundamental fact.
\begin{Theorem} For any principal $\R^\ti$ bundle $M\to M_0$, the tangent $\sT M$ and the co\-tan\-gent~$\sT^* M$ bundle are canonically linear $\R^\ti$-bundles with the $\R^\ti$ action described in Proposition~{\rm \ref{prop:tangent lift}}.
\end{Theorem}
It is easy to see that, starting with coordinates $(t,x^a)$ in $M$, where~$t\in\R^\ti$, identif\/ication~(\ref{ident}) takes in the above cases the form
\begin{gather*}\sT M=M\ti_{M_0}(\sT M/\R^\ti),\qquad \big(t, x^{a}, \dot{t}, \dot{x}^{b}\big)\mapsto \big(t, x^{a}, \dot{\mathbf{t}}, \dot{x}^{b}\big),
\end{gather*}
where $\dot{\mathbf{t}}=t^{-1} \dot t$, and
\begin{gather*}\sT^*M=M\ti_{M_0}\sJ^1M^+,\qquad \big(t,x^a,z,p_b\big)\mapsto \big(t,x^a,z,\mathbf{p}_b\big),\end{gather*}
where ${\mathbf{p}_b}=t^{-1} p_b$.

All this can be easily generalised to a concept of a \emph{graded $\Rt$-bundle}, i.e., a graded bundle \cite{Bruce:2014,Grabowski2012} equipped with a compatible $\Rt$-principal structure. We simply assume that the principal action commutes with the grading represented by a homogeneity structure $\rmh\colon \R\ti P\to P$ (reducing to homotheties in the case of a vector bundle). Also a concept of an \emph{$n$-tuple linear $\Rt$-bundle} is completely obvious: we require the compatibility of the $\Rt$-action with all compatible~$n$ vector bundle structures present on~$P$. For the basics on graded bundles and their relation to $n$-tuple vector bundles we refer to \cite{Bruce:2014,Bruce:2016a,Grabowski2012}.

\subsection[$\Rt$-algebroids vs $\R^\ti$-groupoids]{$\boldsymbol{\Rt}$-algebroids vs $\boldsymbol{\R^\ti}$-groupoids}
Assuming that a Kirillov manifold is equipped simultaneously with a compatible vector bundle, we get the following.
\begin{Definition}[\cite{Grabowski2013}]
A \emph{linear Kirillov structure} is a linear Poisson $\mathbb{R}^{\times}$-bundle, i.e., a linear $\R^\ti$-bundle equipped with a Poisson structure which is linear and homogeneous of degree~$-1$ with respect to the $\R^\ti$-action. If the principal $\mathbb{R}^{\times}$-bundle
is trivial, then we speak about a \emph{linear Jacobi structure}. A~\emph{morphism of linear Kirillov structures} is Poisson morphism that intertwines the respective pairs of~$\R$- and $\R^\ti$-actions.
\end{Definition}

We will denote a linear Kirillov structure as the quadruple $(P, {\rmh}, {\rml},\Lambda)$, where $\rmh$ and $\rml$ are~$\R^\ti$- and~$\R$-actions, respectively, or simply~$(P, \Lambda)$ where no risk of confusion can occur. In local homogeneous coordinates, the Poisson structure must be of the form
\begin{gather*}
\Lambda = \frac{1}{t} \Lambda^{i\alpha}(x) \frac{\partial}{\partial x^{\alpha}} \wedge \frac{\partial}{\partial y^{i}} + \frac{1}{ 2t} y^{k}\Lambda_{k}^{ij}(x) \frac{\partial}{\partial y^{j}} \wedge \frac{\partial }{\partial y^{i}} + \Lambda^{i}(x)\frac{\partial}{\partial y^{i}} \wedge \frac{\partial}{\partial t},
\end{gather*}
where $(t,x,y)$ are coordinates of $(\rmh,\rml)$-bidegrees $(1,0)$, $(0,0)$, and $(0,1)$, respectively.
\begin{Example}
There is a canonical linear Kirillov structure associated with a given Kirillov manifold $(P,\rmh,\zL)$.
It is simply the linear $\Rt$-bundle $\sT P$ (with the tangent lift of $\Rt$-action) equipped with the tangent lift $\dt\zL$ of the Poisson structure $\zL$.
\end{Example}
\begin{Remark}
The above construction is in principle equivalent to the one described in \cite[Remark 2]{Grabowski:2001} and, for trivial $\Rt$-bundles, it leads to the construction of a Lie algebroid associated with a given Jacobi structure, as presented in \cite{Kerbrat:1993}. The above description, however, is strikingly simple.
\end{Remark}

We must draw attention to the similarities with Lie algebroids. In particular Lie algebroid structures on a vector bundle are equivalent to linear Poisson structures on the dual vector bundle; there is an equivalence of categories here. A similar correspondence holds in the case of linear Kirillov structures. One has to take a little care here, as the above picture is dual to the description in terms of homogeneous Poisson structures. In particular, in the above proposition we make the identif\/ication $E = P^{*}$ as vector bundles over $P_{0} = P\slash \mathbb{R}^{\times}$.
The dual vector bundle~$P^*$ is clearly a Lie algebroid which comes with an $\Rt$-action. From considerations preceding \cite[Theorem~8.1]{Grabowski2013} (see also \cite[Theorem~8.2]{Grabowski2013}) it follows easily:
\begin{Proposition}\label{prop:Kirillov Lie}
The objects dual to linear Kirillov structures are precisely $\Rt$-algebroids in the sense of Definition~{\rm \ref{d1}}.
\end{Proposition}

\begin{Remark}
$\Rt$-algebroids associated with trivial principal $\Rt$-structures we call \emph{Jacobi algebroids}. This concept of a Jacobi algebroids coincides with the one introduced and studied in \cite{Grabowski:2001,Grabowski:2003} and that of a \emph{generalised Lie algebroid} in~\cite{Iglesias-Ponte:2001}. In full generality, the notion of a $\Rt$-algebroid is equivalent to that of an \emph{abstract Jacobi algebroid}, def\/ined in \cite{Le:2014} as a Lie algebroid together with a representation thereof on a line bundle, see~\cite{Grabowski2013} for a closer description and proof of this equivalence.
\end{Remark}

\begin{Definition}
An $\Rt$-algebroid is said to be \emph{integrable} if it is integrable as a Lie algebroid.
\end{Definition}

In view of Proposition \ref{prop:Kirillov Lie}, Theorem \ref{thm:lie theory01} immediately implies the following.
\begin{Theorem}\label{thm:lie theory1} There is a one-to-one correspondence between integrable $\Rt$-algebroids and source simply-connected $\mathbb{R}^{\times}$-groupoids.
\end{Theorem}

A generalization of all these concepts to \emph{weighted $\Rt$-algebroids}, i.e., weighted Lie algebroids~\cite{Bruce:2014} with a principal $\Rt$-action by Lie algebroid automorphism commuting with the gradation (homogeneous structure) is straightforward.

\section{Kirillov and contact groupoids}\label{sec:contact and jacobi}
\subsection{Kirillov and Jacobi groupoids}
\begin{Definition}A \emph{Kirillov groupoid} is a $\mathbb{R}^{\times}$-groupoid equipped with a homogeneous multiplicative Poisson structure of degree $-1$, i.e., an $\Rt$-groupoid which has a Poisson groupoid structure of degree $-1$. Kirillov groupoids with trivial $\R^\ti$-bundle will be called \emph{Jacobi groupoids}. If the Poisson structure is non-degenerate, i.e., a symplectic structure, then we will speak of a \emph{contact groupoid}.
\end{Definition}

In a slightly dif\/ferent words, a Kirillov groupoid is a $\mathbb{R}^{\times}$-groupoid $(\mathcal{G}, \rmh)$ equipped with a~multiplicative Poisson structure $\zL$, such that $(\mathcal{G}, \rmh , \zL)$ is a Kirillov manifold (by forgetting the groupoid structure).

\begin{Remark} A contact groupoid is a \emph{homogeneous symplectic groupoid}, i.e., a symplectic groupoid $(\cG,\zw)$ equipped additionally with a compatible principal $\Rt$-bundle structure~$\rmh$ such that $\Rt$ acts by groupoid isomorphisms and $\zw$ is homogeneous of degree $1$ with respect to this action, $\rmh_t^*\zw=t\,\zw$. Symplectic groupoids have been def\/ined by Weinstein \cite{Weinstein:1987} and, under dif\/ferent names, independently by Karasev \cite{Karasev:1987} and Zakrzewski \cite{Zakrzewski:1990,Zakrzewski:1990a}. They can be understood as groupoids $\cG\rightrightarrows M$ equipped with a \emph{multiplicative} symplectic form~$\zw$. The notion of a~homogeneous symplectic groupoid can be traced back to Libermann~\cite{Libermann:1995}, however her notion of homogeneity is in terms of a vector f\/ield and not an action of $\R^\ti$, so does not cover the case of an arbitrary line bundle.
\end{Remark}
\begin{Example} Let $\cG$ be a Lie groupoid. Then, the cotangent bundle $\sT^*\cG$ is canonically a~symplectic groupoid \cite{Weinstein:1987} with respect to the canonical symplectic form $\zw_\cG$ on $\sT^*\cG$. The manifold of units is the dual $\rmA^*(\cG)$ of the Lie algebroid $\rmA(\cG)$ of $\cG$, embedded into $\sT^*\cG$ as the conormal bundle $\zn^*M$. We will refer to it as to the \emph{canonical symplectic groupoid of $\cG$}. It has a vector bundle structure compatible with the groupoid structure in the sense that homotheties \smash{$\rml_t(\zvy_y)=t.\zvy_y$} in the vector bundle $\sT^*\cG\to\cG$ act as groupoid morphisms (it is a \emph{$\mathcal{VB}$-groupoid}). The source and the target maps $s,t\colon \sT^*\cG\to\rmA^*(\cG)$ intertwine the homotheties in $\sT^*\cG\to\cG$ with that in $\rmA^*(\cG)\to M$. It is now clear that removing the level sets $Z_s=s^{-1}(\{ 0\})$ and $Z_t=t^{-1}(\{ 0\})$ gives us an open-dense subgroupoid
\begin{gather*}
\mathcal{C}(\cG)=\sT^*\cG\setminus\{ Z_s\cup Z_t\}\rightrightarrows \rmA^*(\cG)\setminus \{0\}
\end{gather*}
of $\sT^*\cG$. In other words, $\mathcal{C}(\cG)$ consists of covectors from $\sT^*\cG$ which vanish on vectors tangent to source or target f\/ibres. Of course, being an open subgroupoid of $\sT^*\cG$ it is still a symplectic groupoid, but as the zero section of $\sT^*\cG$ has been removed and as $\mathcal{C}(\cG)$ remains $\Rt$-invariant, the group $\Rt$ acts on $\mathcal{C}(\cG)$ by non-zero homotheties in a free and proper way. The symplectic form remains homogeneous of degree $1$ with respect to this action, so we are dealing with a~contact groupoid. The contact groupoid $\mathcal{C}(\cG)$ is canonically associated with the groupoid $\cG$ and will be called the \emph{canonical contact groupoid of $\cG$}. In the traditional picture, it should be viewed as the reduced groupoid $\mathcal{C}(\cG)/\Rt$ which is an open-dense submanifold of the projectivisation bundle~$\mathcal{P}(\sT^*\cG)$.
\end{Example}
\begin{Remark}We will show that Jacobi groupoids in our sense coincide with the Jacobi groupoids def\/ined by Iglesias-Ponte and Marrero in~\cite{Iglesias-Ponte:2003} (also see~\cite{Iglesias-Ponte:2004}), while our contact groupoids are contact groupoids in the sense of Dazord~\cite{Dazord:1995}. The latter are more general and at the same time conceptually simpler than that of Kerbrat and Souici-Benhammadi~\cite{Kerbrat:1993}, which require a globally def\/ined contact form (whose multiplicativity is twisted by a multiplicative function).
\end{Remark}

Since a Kirillov groupoid is both an $\R^\ti$-groupoid and a Poisson groupoid as def\/ined by Weinstein \cite{Weistein:1988}, let us decipher the above def\/inition of the Kirillov groupoid $(\mathcal{G} \rightrightarrows M, {\rmh}, \Lambda)$:
\begin{center}
\leavevmode
\begin{xy}
(0,20)*+{\mathcal{G}}="a"; (20,20)*+{\mathcal{G}_{0}}="b";%
(0,0)*+{M}="c"; (20,0)*+{M_0.}="d";%
{\ar "a";"b"}?*!/_2mm/{\pi};
{\ar@<1.ex>"a";"c"} ;
{\ar@<-1.ex> "a";"c"} ?*!/^3mm/{s} ?*!/_6mm/{t};
{\ar@<1.ex>"b";"d"};%
{\ar@<-1.ex> "b";"d"}?*!/^3mm/{{\zs}} ?*!/_6mm/{{\zt}}; %
{\ar "c";"d"}?*!/^3mm/{\rmp};
\end{xy}
\end{center}
\begin{enumerate}\itemsep=0pt
\item $\zp$ and $\rmp$ are principal $\R^\ti$-bundles.
\item $M$ is a coisotropic submanifold of $(\mathcal{G}, \Lambda)$.
\item There is a unique homogeneous Poisson structure on $M$ such that the source map is a Poisson morphism and the target map an anti-Poisson morphism. Thus $(M, \rmh ,s_{*}\Lambda )$, where $\rmh$ is the homogeneity structure restricted to $M$, is a Kirillov manifold.
\item The Poisson structure induces a morphism of Lie groupoids:
\begin{center}
\leavevmode
\begin{xy}
(0,20)*+{\sT^{*}\mathcal{G}}="a"; (20,20)*+{ \sT \mathcal{G}}="b";%
(0,0)*+{\rmA^{*}(\mathcal{G})}="c"; (20,0)*+{\sT M.}="d";%
{\ar "a";"b"}?*!/_3mm/{\Lambda^{\#}};
{\ar@<1.ex>"a";"c"};
{\ar@<-1.ex> "a";"c"};
{\ar@<1.ex>"b";"d"};%
{\ar@<-1.ex> "b";"d"}; %
{\ar "c";"d"}?*!/^3mm/{\rho_\zL};
\end{xy}
\end{center}
According to Theorem \ref{liftingactions}, $\Lambda^{\#}$ intertwines the lifted actions of $\R^\ti$, so induces a morphism of the reduced groupoids
\begin{gather*}
\Lambda^{\#}_0\colon \ \sT^{*}\mathcal{G}/\R^\ti\rightarrow \sT \mathcal{G}/\R^\ti.
\end{gather*}
\item According to Remark \ref{redu}, in the case of the trivial $\R^\ti$-bundle, $\cG=\cG_0\ti\R^\ti$, we have the identif\/ications
 \begin{gather*}\sT^{*}(\mathcal{G}_0\ti\R^\ti)/\R^\ti\simeq \sT^{*}\mathcal{G}_0\ti\R,\qquad
 \sT(\mathcal{G}_0\ti\R^\ti)/\R^\ti\simeq \sT\mathcal{G}_0\ti\R,
 \end{gather*}
 and the reduced morphism can be viewed as a map
\begin{gather*}\Lambda^{\#}_0\colon \ \sT^{*}\mathcal{G}_{0} \times \mathbb{R} \rightarrow \sT \mathcal{G}_{0} \times \mathbb{R},
\end{gather*}
associated with a Jacobi structure on $\cG_0$. This map is a groupoid morphism for
the groupoid structures determined by that on $\cG_0\ti\Rt$, so according to Theorem \ref{trivialsplit}, by the groupoid structure of $\cG_0$ and a~multiplicative function $\rmb\colon \cG_0\to\Rt$, which reduces to a~multiplicative function $\log{|\rmb|}$ into the additive group of reals. This leads to the def\/inition of a~Jacobi groupoid as presented in \cite{Iglesias-Ponte:2003}, although the explicit form of the groupoid structures on $\sT^{*}\mathcal{G}_{0} \times \mathbb{R}$ and $\sT\mathcal{G}_{0} \times \mathbb{R}$, expressed in terms of $\cG_0$ and $\rmb$, is quite complicated. Our `Kirillov version' of the Jacobi groupoid is not only more general, but conceptually simpler. The technical complications of the def\/inition in \cite{Iglesias-Ponte:2003}, together with the presence of the multiplicative function $\log|\rmb|$, come from insisting on working with trivialisations of principal bundles.
\end{enumerate}

The above observations do not change if we consider a contact groupoid; in particular we still have a Kirillov manifold $M$ which is generally not a contact manifold.

\begin{Definition}
A Kirillov manifold is said to be an \emph{integrable Kirillov manifold} if it arises from a contact groupoid as described above.
\end{Definition}

By minor modif\/ication of the classical results on the integrability of Poisson manifolds and in view of Theorem~\ref{thm:lie theory1}, we are led to the following;
\begin{Proposition}\label{prop:lie theory2}
The following statements are equivalent:
\begin{enumerate}\itemsep=0pt
\item[$1.$] The Kirillov manifold $(M, \rmh, \Lambda)$ is integrable.
\item[$2.$] The $\Rt$-algebroid $(\sT M, \sT \rmh, \rmd_{\sT}\Lambda)$ is integrable.
\item[$3.$] The Lie algebroid structure on $\sT^*M$ corresponding to the linear Poisson structure $\dt\zL$ is integrable.
\item[$4.$] The Poisson structure $\Lambda$ is integrable.
\end{enumerate}
\end{Proposition}

\begin{Remark} Kerbrat and Souici-Benhammadi \cite{Kerbrat:1993}, \emph{define} the integrability of a Jacobi mani\-fold $\big(M , \widetilde{\zL}, \widehat{R}\big)$ as the integrability of the associated Lie algebroid $\sT^*M \times \R $. We will not give a careful def\/inition here of this Lie algebroid structure and simply point the reader to the ori\-gi\-nal literature. One of the main theorems of Crainic and Zhu \cite[Theorem~1(iii)]{Crainic:2007} is that the integrability of a Jacobi manifold is equivalent (though non-trivially) to the integrability of its poissonisation. From our point of view, it is the integrability of the Kirillov manifold associated with a Jacobi manifold, i.e., its poissonisation, that is fundamental and the `correct' starting def\/inition of the integrability of a Jacobi manifold.
\end{Remark}

\subsection{Contact groupoids}
Let us restrict attention to a contact groupoid $(\cG,\rmh,\zw)$. We know, \emph{via} Remark \ref{rem:contact distribution}, that the homogeneous symplectic structure $\zw$ on the $\Rt$-bundle $(\cG,\rmh)$ is equivalent to a contact structure $C=C(\cG,\rmh,\zw) \subset \sT^{*}\mathcal{G}_{0}$, and further to a contact distribution $\mathcal{D}=\mathcal{D}(\cG,\rmh,\zw)=C^0 \subset \sT \mathcal{G}_{0}$. Such a distribution is called a \emph{contact groupoid} by Dazord \cite{Dazord:1995,Dazord:1997} (and \emph{conformal contact groupoid} by the authors understanding contact groupoids as groupoids equipped with a globally def\/ined contact form) if the contact distribution is closed with respect to the operation in the tangent groupoid $\sT\cG_0$: it is invariant with respect to inversion and $\mathcal{D}\bullet\mathcal{D}\subset\mathcal{D}$, where $\bullet$ is the (partial) multiplication in $\sT\cG_0$. In other words, $\mathcal{D}$
is a Lie subgroupoid of $\sT\cG_{0} \rightrightarrows \sT M$. We refer to such Lie groupoids as \emph{Dazord groupoids}, keeping the term contact groupoids to refer to our notion. However, the two notions are equivalent.
\begin{Theorem}\label{thm:Dazord}
Any contact groupoid $(\cG,\rmh,\zw)$, with $\cG/\Rt=\cG_0$, has a canonical and equivalent realisation as each of the following:
\begin{itemize}\itemsep=0pt
\item a contact subgroupoid $C^\ti(\cG,\rmh,\zw)$ of the canonical contact groupoid $\mathcal{C}(\cG_0)$;
\item a Dazord groupoid $\mathcal{D}(\cG,\rmh,\zw)$, being simultaneously a contact distribution and a subgroupoid of $\sT\cG_0$.
\end{itemize}
\end{Theorem}
\begin{proof} Let $\nabla$ be the Euler vector f\/ield on $\cG$. Since $\nabla$ generates groupoid isomorphisms, it is a~mul\-tiplicative vector f\/ield on $\cG$, thus $\nabla\colon \cG\!\to\!\sT\cG$ is a groupoid morphism (over \smash{$\nabla_{|M}\colon M\!\to\!\sT M$}). Since the symplectic form $\zw$ is multiplicative, it def\/ines an isomorphism of groupoids $\zw^\flat\colon \sT\cG\to\sT^*\cG$. The one-form $\zh=\zw^\flat\circ\nabla\colon \cG\to\sT^*\cG$ is a groupoid morphism, thus multiplicative, $\zh(yy')=\zh(y)\star\zh(y')$, where $\star$ is the groupoid multiplication in the cotangent groupoid $\sT^* \cG\rightrightarrows \rmA^*(\cG)$. We know that it takes values in basic covectors, $\zh(y)=\zp^*(\Psi(y))\in\sT^*_y\cG$, $\Psi(y)\in\sT^*_{\zp(y)}\cG_0$, so can be viewed as a map $\Psi\colon \cG\to\sT^*\cG_0$. Consequently, the range $C^\ti(\cG,\rmh,\zw) =\{\Psi(y)\colon y\in\cG\}$ of $\Psi$ is a~subgroupoid in the cotangent groupoid $\sT^*\cG_0$. According to Theorem~\ref{thm:contact}, $\Psi$ is also an embedding of contact structures (homogeneous symplectic $\Rt$-bundles), so $\Psi$ is just a realisation of $(\cG,\rmh,\zw)$ as a contact subgroupoid of $\mathcal{C}(\cG_0)$.

Note f\/inally that the contact distribution $\mathcal{D}(\cG,\rmh,\zw)$ and the contact structure $C(\cG,\rmh,\zw)$ are related by the polar condition: one annihilates the other in the canonical pairing between the tangent and the cotangent bundle. Since the partial multiplication in the cotangent and tangent groupoid are related by the condition
\begin{gather*}
\zvy_g\star\zvy'_h(X_g\bullet X'_h)=\zvy_g(X_g)+\zvy'_h(X'_h) ,\end{gather*}
it can be easily seen that $C(\cG,\rmh,\zw)$ is a subgroupoid if and only if $\mathcal{D}(\cG,\rmh,\zw)$ is a subgroupoid.
\end{proof}

\begin{Remark} To be very clear, $\Psi\colon \mathcal{G} \rightarrow \sT^{*}\mathcal{G}_{0}$ is in general \emph{not} a~Lie groupoid morphism (a~contact form on $\cG_0$ need not to be multiplicative), however the range of $\Psi$ is a contact subgroupoid of $\sT^{*}\mathcal{G}_{0}$. Thus we have a canonical realisation of $\mathcal{G}$ rather than a genuine morphism between contact groupoids.
\end{Remark}

\subsection*{Acknowledgements}

The authors are indebted to the anonymous referees whose comments have served to improve the content and presentation of this paper. The research of K.~Grabowska and J.~Grabowski was funded by the Polish National Science Centre grant under the contract number DEC-2012/06/A/ST1/00256.

\pdfbookmark[1]{References}{ref}
\LastPageEnding

\end{document}